\documentclass{amsart}

\usepackage{mathabx,hyperref,mathtools,mathrsfs,amsrefs,lineno}
\usepackage{amsfonts,amssymb}
\usepackage{color,url,hyperref}
\usepackage{tikz}
\usetikzlibrary{cd,arrows}

\DeclareSymbolFont{bbold}{U}{bbold}{m}{n}
\DeclareSymbolFontAlphabet{\mathbbold}{bbold}

\newcommand\N{{\mathbb{N}}}
\newcommand\Z{{\mathbb{Z}}}

\newtheorem{thm}{Theorem}[section]

\theoremstyle{definition}

\theoremstyle{remark}

\begin{document}
\author{Laurent Bartholdi} 
\email{laurent.bartholdi@gmail.com}

\author{Roman Mikhailov}

\title{The topology of poker}

\begin{abstract}
  We examine the complexity of the ``Texas Hold'em'' variant of poker from a topological perspective. We show that there exists a natural simplicial complex governing the multi-way winning probabilities between various hands, and that this simplicial complex has dimension at least $4$. We deduce that evaluating the strength of a pair of cards in Texas Hold'em is an intricate problem, and that even the notion of who is bluffing against whom is ill-defined in some situations. The use of topological methods to study intransitivity of multi-player games seems new.
\end{abstract}



\thanks{L.B. is partially supported by ERC AdG grant 101097307}

\date{\today}

\maketitle

\section{Introduction}
In the popular ``Texas Hold'em'' variant of poker (see e.g.~\cite{MR2656351}*{Chapter~22}), you and each of your opponents are dealt two cards, and five cards will be dealt to the table. The winner is the player making the best $5$-card poker game out of their's and the table's cards. Suppose you hold $J\clubsuit10\clubsuit$ and two other players respectively hold $2\diamondsuit2\heartsuit$ and $K\clubsuit2\clubsuit$. Who is favourite? And what happens after one of the opponents folds?

Knowing the winning probabilities of a hand against another one is fundamental to any poker strategy, and are at the heart of von Neumann's analysis of poker~\cite{zbMATH02577401}. What we argue, however, is that winning probabilities \emph{give at best partial information on the current game state}, and sometimes \emph{are meaningless}.

Let us pause to consider a much simpler game, ``Rock, Paper, Scissors'' (RPS). It would be absurd, in a televised retransmission of a RPS match, to display winning probabilities for each player, since by the game's symmetry each player wins against one play and loses against another. The situation, in poker, is \emph{far worse}, due both to the richness of the game and to the high number of players (typically, $8$).

But first, the situation of RPS \emph{does} occur in poker: for the pairs $J\clubsuit10\clubsuit$, $2\diamondsuit2\heartsuit$ and $K\clubsuit2\clubsuit$ we can check that the $1$-on-$1$ winning chances are
\[w(J\clubsuit10\clubsuit,2\diamondsuit2\heartsuit)=54\%,\quad w(2\diamondsuit2\heartsuit,K\clubsuit2\clubsuit)=63\%,\quad w(K\clubsuit2\clubsuit,J\clubsuit10\clubsuit)=55\%.\]
This means that the \emph{same player} may become favourite or underdog depending on which of its opponents folds.

To model such seemingly paradoxical phenomena, we introduce a topological invariant of a game, and apply computational tools to derive the first non-trivial results for poker. We avoid any discussion on the exact winning probabilities to concentrate only on \emph{who} wins: we therefore have a set $X$ of player hands, and an \emph{antisymmetric} relation $r\subseteq X\times X$, namely a relation such that for all $x,y\in X$ at most one of $r(x,y),x=y,r(y,x)$ holds. An antisymmetric relation $r\subseteq X\times X$ is a partial order if and only if it is \emph{transitive}: $X$ does not contain elements $R,P,S$ with $r(R,P),r(P,S),r(S,R)$.

This property can be interpreted topologically, by means of a \emph{simplicial complex}, namely a collection of subsets of $X$ that is closed under taking subsets. This is the abstract formulation of a union of copies of standard $n$-simplices $\{x\in[0,1]^{n+1}:\sum x_i=1\}$, and indeed every simplicial complex $\mathscr K_X$ on $X$ has a geometric realization as $\{x\in[0,1]^X\mid\text{support}(x)\in\mathscr K_X,\sum x_i=1\}$.

Consider the simplicial complex $\mathscr K_X$ with vertex set $X$ and simplices the ordered subsets of $(X,r)$, namely $C\subseteq X$ is a simplex if and only if the restriction of $r$ to $C\times C$ is transitive. This construction is classical when starting with a partial order --- the theory of simplicial complexes and of partially ordered sets are essentially equivalent, and see~\cite{MR1373690} for its applications to topology --- but nothing prevents us from applying it to our relation $r$, and studying the topology of $\mathscr K_X$. For example, if $r$ is \emph{dominated} in the sense that for all $x,y\in X$ there exists $z$ with $r(x,z)$ and $r(y,z)$, then $\mathscr K_X$ is contractible.

Simple game-theoretical properties may then be translated to topology. For example, if one hand $x\in X$ beats all the others, then $\mathscr K_X$ is a cone with apex $x$; and the simplicial complex associated with ``Rock, Paper, Scissors'' is homeomorphic to the circle $S^1$.

We posit that the complexity of a game is related to the topological complexity of $\mathscr K_X$. As a simple proxy, we consider the maximal homological dimension of a subcomplex; a measure that both captures the topological richness and is computable, or at least boundable from below. By contrast, the Lusternik-Schnirelmann category seems even harder to compute, and less meaningful.

\subsection{Main results}
We consider the set $X=\{(i,j):1\le i<j\le 52\}$ of pairs of cards, and the relation
\[r((i,j),(k,\ell))\Longleftrightarrow\begin{cases}\text{averaging over the $\binom{52-4}{5}$ remaining cards,}\\\text{pair $(i,j)$ wins with more than $50\%$ chance}.\end{cases}\]
We sometimes write `$(i,j)>(k,\ell)$' for $r((i,j),(k,\ell))$, without implying transitivity.

Our main result is that the homotopy type of Texas Hold'm is quite rich; in other words, there are intricate card configurations in which every player could be winning against another one; thus even the concept of ``bluff'' (see e.g.~\cite{MR3399101}) needs to be revisited since it is impossible to define, at some moments, who is bluffing against whom:
\begin{thm}
  For Texas Hold'em, the simplicial complex $\mathscr K_X$ contains $4$-dimensional subcomplexes.
\end{thm}
More precisely, we exhibit $S^4$ as a subcomplex of $\mathscr K_X$.

If there were a $4$-dimensional simplex in $\mathscr K_X$, it could be interpreted as follows: \emph{there is a configuration with $5$ players such that no hand is better than the others, but as soon as a player folds the remaining $4$ are linearly ordered}. Loosely speaking, our result shows that the same phenomenon occurs with coalitions instead of players.

There is a long history of associating numbers to games; integers for Sprague and Grundy~\cites{sprague,grundy} or generalized (and in particular not necessarily ordered) numbers for Conway~\cite{MR1803095}. We view topological invariants as a new powerful tool supplementing these more classical invariants.

\subsection{Acknowledgments}
The calculations have made heavy use of the computer algebra program \texttt{Oscar}~\cite{OSCAR}, as well as the poker hand evaluator~\texttt{PokerHandEvaluator.jl}.

\section{Implementation}
Using the computer language \textsf{Julia} and its packages \textsf{PlayingCards} and \textsf{PokerHandEvaluator}, we computed the relation $r$. (We also independently re-implemented \textsf{PokerHandEvaluator} to make sure of its correctness.) The code is available on the Zenodo repository \url{https://doi.org/10.5281/zenodo.7885276}; the file \texttt{poker-data.jl} defines an array \texttt{CARDPAIRS} of size $1326$ listing all pairs of cards, and the array \texttt{r} in the HDF5 dataset \texttt{poker-data.hdf5} has size $1326\times1326\times3$, in such a way that \texttt{r[i,j,1:3]=(w,t,l)} means that playing \texttt{CARDPAIRS[i]} against \texttt{CARDPAIRS[j]} results in \texttt{w} wins, \texttt{t} ties and \texttt{l} losses when considering all $\binom{48}{5}$ possible table cards; thus $\texttt w+\texttt t+\texttt l=1712304$.

We then explored subsets $Y\subseteq X$, computed the corresponding simplicial complex $\mathscr K_Y$ using the computer algebra package \textsf{Oscar} and its \textsf{Polymake} interface, and its homology.

Let us begin with an example of an $S^1$ in $\mathscr K_X$. Consider the hands $A\clubsuit2\clubsuit$, $3\clubsuit5\clubsuit$, $2\diamondsuit2\heartsuit$. Note that $3\clubsuit5\clubsuit$, known as ``Carabas'' in Russian, is a well-known tricky holding. The first wins against the second on average, because of the strength of the ace. The second wins against the third because of the possibilities of forming a flush. The third wins against the first because of the pair. These winning probabilities are respectively $0.591$, $0.504$, $0.620$. We write these data in the following diagram:
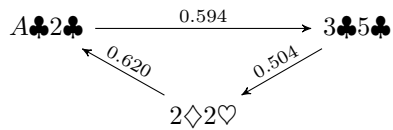
\[\begin{tikzcd}[>=stealth']
    A\clubsuit2\clubsuit \ar[->,rr,"0.594"] & & 3\clubsuit5\clubsuit\\ & 2\diamondsuit2\heartsuit \ar[->,ul,sloped,"0.620"]\ar[<-,ur,below,sloped,"0.504"]
  \end{tikzcd}
\]  

It seemed useful to consider hands in which the cards are close, so we ordered the $52$ cards by their value, kept each card independently with some probability $p$, and formed the pairs out these cards in increasing order. After a few thousand runs, our computer search came up with the subset
\[Y=\{A\clubsuit A\diamondsuit,\; A\heartsuit A\spadesuit,\quad
  6\diamondsuit6\spadesuit,\; J\spadesuit Q\diamondsuit,\; 10\heartsuit J\heartsuit,\quad
  2\heartsuit2\spadesuit,\; 7\spadesuit 10\clubsuit,\; 4\clubsuit6\clubsuit\}.
\]

The relations between these pairs are given in the following diagram:
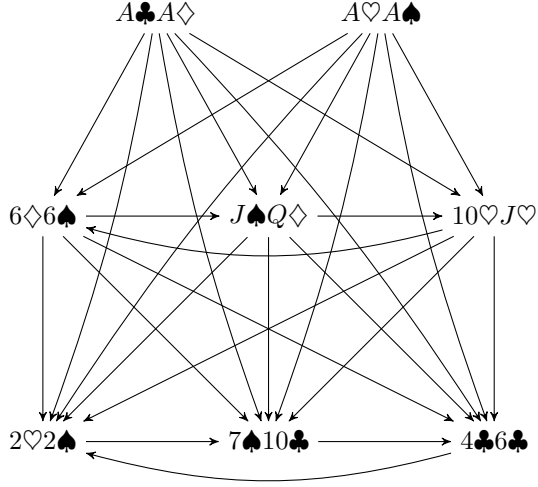
\[\begin{tikzpicture}[>=stealth']
    \node (1) at (-1.5,2.7) {$A\clubsuit A\diamondsuit$};
    \node (2) at (1.5,2.7) {$A\heartsuit A\spadesuit$};
    \node (3) at (-3,0) {$6\diamondsuit6\spadesuit$};
    \node (4) at (0,0) {$J\spadesuit Q\diamondsuit$};
    \node (5) at (3,0) {$10\heartsuit J\heartsuit$};
    \node (6) at (-3,-3) {$2\heartsuit2\spadesuit$};
    \node (7) at (0,-3) {$7\spadesuit 10\clubsuit$};
    \node (8) at (3,-3) {$4\clubsuit6\clubsuit$};
    \draw[->] (1) edge (3) edge (4) edge (5) edge[bend left=5] (6) edge[bend right=5] (7) edge[bend left=5] (8);
    \draw[->] (2) edge (3) edge (4) edge (5) edge[bend right=5] (6) edge[bend left=5] (7) edge[bend right=5] (8);
    \foreach\i in {3,4,5}\foreach\j in {6,7,8} {\draw[->] (\i) edge (\j); }
    \draw[->] (3) edge (4) (4) edge (5) (5) edge[bend left=15] (3);
    \draw[->] (6) edge (7) (7) edge (8) (8) edge[bend left=15] (6);
  \end{tikzpicture}
\]

A direct calculation showed that all its homotopy groups are trivial except $H_4(\mathscr K_Y,\Z)=\Z$.

In fact, this can be checked manually, assuming of course knowledge of the winning relation $r$. Indeed the first row defines an $S^0$, since these card holdings are incomparable; the second and the third row define $S^1$, with $6\diamondsuit6\spadesuit>J\spadesuit Q\diamondsuit>10\heartsuit J\heartsuit>6\diamondsuit6\spadesuit$ and $ 2\heartsuit2\spadesuit>7\spadesuit 10\clubsuit>4\clubsuit6\clubsuit>2\heartsuit2\spadesuit$; and the complex $\mathscr K_Y$ is the join of these three subcomplexes, totally ordered as $\{A\clubsuit A\diamondsuit,A\heartsuit A\spadesuit\}>\{6\diamondsuit6\spadesuit,J\spadesuit Q\diamondsuit,10\heartsuit J\heartsuit\}>\{2\heartsuit2\spadesuit,7\spadesuit 10\clubsuit,4\clubsuit6\clubsuit\}$, hence is homeomorphic to $S^4$.

\section{Outlook}
We have barely scratched the surface of the topological complexity of Texas Hold'em. In particular, it does not seem possible to compute the homotopy type, or even just the homology, of $\mathscr K_X$ with current technology.

Indeed, using the usual limit of $10$ players per table, there are $\binom{52\cdot51/2}{10}$ collections of pairs of hands to consider, and for each a homological calculation to perform. This homological calculation is feasible in the cases we considered, but can become quite expensive, since in general computing homology groups is NP-hard~\cite{MR3504911}. When considering such large datasets, researchers typically concentrate on $H_1$ and possibly $H_2$, while we are interested in higher-degree phenomena.

\subsection{Revealing table cards one at a time}
We have made the simplifying assumption that all table cards are unknown. In actual Texas Hold'em, there are more than one bidding round, and more cards are progressively revealed. As more cards are revealed, $X$ shrinks to a subset $X'$ because fewer cards may appear in our opponents' hands; it would be interesting to study the importance of the partial revealing of information, in the form of a failure of the inclusion $X'\hookrightarrow X$ to induce a simplicial map.

\subsection{Sensitivity to data}
As shown by the brief calculation above, the probabilities associated with edges in $r$ are typically not microscopically away from $0.5$. The closest one is $3x3y$ versus $Ax10x$ for two suits $x\neq y\in\{\clubsuit,\diamondsuit,\heartsuit,\spadesuit\}$, with winning probability $0.50007$. The natural tool with which to explore the sensitivity of $\mathscr K_X$ is \emph{persistent homology}: for a probability $p\in[0.5,1]$, consider the relation $r_p$ in which $r_p((i,j),(k,\ell))$ means that the probability that $(i,j)$ wins is at least $p$, and the associated simplicial complex $\mathscr K_{X,p}$. We have for each $0.5\le p<q\le1$ inclusion maps $\mathscr K_{Y,q}\hookrightarrow\mathscr K_{Y,p}$; what are their relative homologies, as $Y$ and $p,q$ vary?

\subsection{Other intransitive games}
There is a substantial literature on ``intransitive games'', see~\cite{MR1880801}*{Chapters~22 and~23}, however only considering $2$ players. One of the most promising ones is the ``Penney game''~\cite{MR0611250}. Some parameter $n\in\N$ is fixed. Every player chooses a binary sequence of length $n$. An infinite binary sequence is then drawn at random, one bit at a time. The first player whose sequence shows up wins.

It is well known that this game is not transitive; for $n=3$, we have $011>110>100>001>011$, and the associated complex $\mathscr K_{\{0,1\}^3}$ is a bouquet of $3$ circles. We computed the homology of the whole complex $\mathscr K_{\{0,1\}^n}$ for $n\le6$, giving for the last case
\[H_*(\mathscr K_{\{0,1\}^6},\Z)=(0,0,0,0,0,\Z^{38},\Z^{149},\Z^{12})\text{ for }0\le*\le7.\]
There does not seem to be any obvious pattern to these numbers.

Poker is quite apart from these games in that it is a real-life game, with a large population of expert or professional players, in which it is essential to estimate with accuracy the winning chances relative to all the other participants. Our results show that these data must be considered globally, and that the one-on-one probabilities only serve as the carrier for a powerful, high-dimensional topological invariant.


\end{document}